\documentstyle[12pt]{article}
\newcommand{\N}{\mbox{I$\!$N}}
\newcommand{\R}{\mbox{I$\!$R}}
\newcommand{\Hb}{\mbox{l$\!$H}}

\newcommand{\qed}{{\hfill
{$\rlap{$\sqcap$}\sqcup$}}\\[0.2in]\hspace*{0.5in}}
\newcommand{\qedwh}{{\hfill {$\rlap{$\sqcap$}\sqcup$}}\\[0.2in]}
\newcommand{\bk}{\\[0.03in] \hspace*{0.5in} }
\newcommand{\btd}{\bigtriangledown}

\begin{document}

\vspace*{0.07in}
\begin{center} {\LARGE   Uniqueness of Positive Solutions of the }
\medskip \\ {\LARGE  Equation $\Delta_g u + c_n u = c_n u^{{n + 2}\over {n
- 2}}$ \ and}
\medskip \\ {\LARGE Applications to Conformal Transformations } 
\medskip \medskip \medskip \medskip\\   {\large Man Chun LEUNG}
\medskip\\
{\large Department of Mathematics, National University of Singapore,} \\
{\large
Singapore 119260  \ \ {\tt matlmc@nus.sg} \ \ Fax. No. 65-779-5452}
\medskip \medskip \\ 
\end{center}

\begin{abstract}We study uniqueness of positive solutions to the
equation
$\Delta_g u + c_n u = c_n u^{{n + 2}\over {n - 2}}$ on complete Riemannian
manifolds. We apply the results to show
that conformal transformations on certain complete Riemannian manifolds of
constant negative scalar curvature are isometries. We also study
uniqueness of
complete positive solutions and radial solutions.
\end{abstract}
\vspace{0.2in} 
KEY WORDS \& AMS MS Classifications\footnote{
\hspace*{0.24in}1991 AMS MS Classifications: Primary 35J60; Secondary
58G03
53C20}

\vspace{0.5in}

{\bf \Large {\bf 1. \ \  Introduction}}

\vspace{0.3in}

In this paper we study uniqueness of positive solutions to the equation
$$\Delta_g u + c_n u = c_n u^{{n + 2}\over {n - 2}}\,, \leqno (1.1)$$ 
where $c_n = (n - 2)/[4 (n - 1)]$ and $\Delta_g$ is the Laplacian operator
for  
a Riemannian manifold $(M, g)\,.$ Equation (1.1) arises as in the
conformal
deformation of a Riemannian metric into constant negative scalar
curvature. Let
$M$ be an open
$n$-manifold with
$n
\ge 3\,.$ Bland and Kalka
\cite{Bland-Kalka} show that there exists a complete Riemannian metric $g$
on
$M$ with scalar curvature equal to $-1\,.$ Let $u$ be a positive smooth
function
on
$M\,.$ If the conformal metric
$g_c = u^{4\over {n - 2}} g$ has constant scalar curvature equal to
$-1\,,$ then
$u$ satisfies equation (1.1). In this paper we study when a positive
solution to
equation (1.1) is unique, that is, whether $u \equiv 1$ on $M\,.$\bk 
Let $(M, g_o)$ be a complete non-compact Riemannian $n$-manifold with
scalar
curvature $R_{g_o}.$ Aviles and McOwen \cite{A-M} study when there exists
a
positive smooth function $v$ on $M$ such that the conformal metric $g =
v^{4\over {n - 2}} g_o$ has scalar curvature equal to $-1\,.$ In this case
$v$
satisfies the equation 
$$\Delta_{g_o} v - c_n R_{g_o} v = c_n v^{{n + 2}\over {n - 2}}\,, \leqno
(1.2)$$
where $R_{g_o}$ is the scalar curvature of $(M, g_o)\,.$ The question
whether
a positive solution to equation (1.2) is unique is equivalent to the
uniqueness
of a positive solution to equation (1.1). For if
$w$ is a smooth function on
$M$ which satisfies equation (1.2), then $u =w/v$ satisfies the
equation (1.1). Therefore $u \equiv 1$ if and only if $v \equiv w$ on
$M\,.$\bk
Equation (1.1) is a semilinear elliptic equation involving the critical
Sobolev
exponent. For open domains in ${\R}^n\,,$ semilinear elliptic equations
which involve critical Sobolev exponents are studied by Brezis and
Nirenberg
\cite{Brezis-Nirenberg},   Gidas, Ni and Nirenberg
\cite{Gidas-Ni-Nirenberg},
Peletier and Serrin
\cite{Peletier-Serrin}, Ni and Serrin
\cite{Ni-Serrin} and Zhang \cite{Zhang.1} \cite{Zhang.2}. In case
$M$ is a compact manifold without boundary, the maximum principle implies
that
equation (1.1) has only one positive solution $u
\equiv 1\,.$ If $M$ is  an open manifold, then equation (1.1) may have
more than
one positive solution. An example can be given on the open unit ball on
${\R}^n$
with the hyperbolic metric. However if the solution is strictly positive
we have
the following result.\\[0.2in]
{\bf Theorem A.} \ \ {\it For an integer $n \ge 3\,,$ let $(M, g)$ be a
complete
Riemannian $n$-manifold with scalar curvature equal to $-1\,.$ Assume that
the
sectional curvature of $(M, g)$ is bounded from below and the injectivity
radius
of $(M, g)$ is positive. If $u \in C^\infty (M)$ is a positive solution to
equation (1.1) with $\inf_{x \in M} u (x) > 0\,,$ then $u \equiv 1$ on
$M$.}\\[0.2in]
\hspace*{0.5in}Theorem A has nice applications to conformal
transformations. Let
$$F : (M, g) \to (M, g)$$  
be a conformal transformation, where the scalar
curvature of $(M, g)$ is equal to $-1\,.$ Let
$u = |F'|^{{n - 2}\over 2}\,,$ where $|F'|$ is the linear stretch factor.
Then
$u$ satisfies equation (1.1) \cite{Schoen}. Therefore uniqueness of 
positive solutions to (1.1) implies that $F$ is an isometry. Using theorem
A we
can prove the following result.\\[0.2in]
{\bf Theorem B.} \ \ {\it For an integer $n \ge 3\,,$ let
$(M, g)$ be a complete non-compact Riemannian $n$-manifold with constant
negative
scalar curvature. Assume that the sectional curvature of
$(M, g)$ is bounded from below and the injectivity radius of $(M, g)$ is
positive. Then any conformal transformation
$F: (M, g) \to (M, g)$ is an isometry.}\\[0.2in]
\hspace*{0.5in}Kuiper \cite{Kuiper} proves that any conformal
transformation on a
complete Einstein manifold with negative scalar curvature is an isometry. 
Lafontaine \cite{Lafontaine} shows that any conformal transformation of
$S^n
\times {\Hb}^m$ ($m
\ge 2$) is an isometry. As a result of theorem B, we can show that any
conformal
transformation of a locally symmetric space of non-compact type with
dimension
bigger than two is an isometry.\bk
The assumption that $u$ is bounded away from zero may not be natural in
the
general content. But there may have more than one positive solution to
equation
(1.1) if we take away the assumption. For example, for any constant $b \ge
1\,,$
$$u_b (x) = \left( {{b \,(1 - |x|^2)}\over { b^2 - |x|^2}}
\right)^{{n - 2}\over 2}  \ \ \ \ {\mbox{for}} \ \ \ x \in B^n $$ 
is a positive solution to equation (1.1) on the (scaled) hyperbolic space.
A more
natural assumption is that the conformal metric $g_c = u^{4\over {n - 2}}
g$ is
{\it complete}. We study the uniqueness of a positive radial solutions to
equation (1.1) under this assumption. For an integer $n \ge 3\,,$ let $(r,
\Theta)$ be the polar coordinates on ${\R}^n\,,$ where $r \ge 0$ and
$\Theta \in
S^{n - 1}\,.$ We consider warped product metrics on ${\R}^+ \times S^{n -
1}$. 
Let $f$ be a smooth function on
$[0, \infty)$ with
$f (0) = 0\,, f' (0) = 1$ and $f > 0$ on
$(0, \infty)\,.$ The Riemannian metric 
$$g = dr^2 + f^2 (r) \,d\Theta^2 \leqno (1.3)$$
is complete on ${\R}^n \cong {\R}^+ \times S^{n - 1}$. The hyperbolic
space is
corresponding to 
$$f (r) = \sinh r \ \ \ \ {\mbox{for}} \ \ r \ge 0\,,$$ 
while the Euclidean space
with the standard metric is corresponding to $f (r) = r$ for all
$r
\ge 0\,.$ Let $u$ be a solution to equation (1.1) with $g$ of the form
(1.3)\,. We say that $u$ is a radial solution if $u$ is a function of $r$
only (independent on $\Theta$).
\\[0.2in]   {\bf Theorem C.} \ \ {\it
Let $u$ be a positive radial
$C^2$-solution to  equation (1.1) on}
$$({\R}^n\,, \ g = dr^2 + f^2 (r) \,d\Theta^2)\,.$$ 
{\it Assume that $\lim_{r \to \infty} f (r) = \infty$ and there
exist positive constants
$R_o$ and $C$ such that 
$|f' (r)/ f (r)| \le C$ for all
$r > R_o\,.$ If the conformal metric
$g_c = u^{4\over {n - 2}} g$ is complete, then
$u \equiv 1\,.$}\\[0.2in]
\hspace*{0.5in}For a positive solution $v$ to equation (1.1)
that may not be radial, we show that if there exists a $\delta > 0$ such
that $v
\le 1 -
\delta$ outside a compact set, then the conformal metric $g_c = v^{4\over
{n -
2}} g$ is not complete. We also show that if $u$ is a positive solution to
equation (1.1) on the (scaled) hyperbolic space and $u \not\equiv 1\,,$
then
$u < 1$ everywhere. 

\pagebreak

{\bf \Large {\bf 2. \ \  Bounded Solutions}}

\vspace{0.3in}

Throughout this paper we let $(M, g)$ be a complete Riemannian
$n$-manifold with
$n \ge 3\,.$ Let ${\mbox{Ric}}_g$ be the Ricci
curvature of
$(M, g)\,.$\\[0.2in] {\bf Lemma 2.1.} \ \ {\it Assume that ${\mbox{Ric}}_g
\ge -
c^2 g$ on $M$ for some positive constant $c\,.$ Let $u \in C^\infty (M)$
be a
positive solution to (1.1). If $\,\inf_{x \in M} u (x) > 0$ and $\,\sup_{x
\in
M} u (x) < \infty\,,$ then
$u \equiv 1$ on $M\,.$}\\[0.1in]
{\it Proof.} \ \ Suppose that there exists a point $x_o \in M$ such that
$u
(x_o) > 1\,.$ Then we have 
$$1 < \sup_{x \in M} u (x) < \infty\,.$$
By the Omori-Yau maximum principle (c.f. \cite{Ratto-Rigoli-Setti}), there
exists
a sequence
$\{ x_k
\}_{k
\in {\N}}$ such that 
$$\lim_{k \to \infty} u (x_k) = \sup_{x \in M} u (x)\,,  \ \ \ \ 
|\btd u (x_k)|
\le {1\over k}$$ 
and 
$$\Delta_g u (x_k) \le {1\over k}\,. \leqno (2.2)$$
As $\sup_{x \in M} u (x) > 1\,,$ there exists a number $\epsilon > 0$ and
$N \in
{\N}$ such that 
$$ c_n \left( u^{{n + 2}\over {n - 2}} (x_k) - u (x_k) \right) > \epsilon
\ \ \
\ {\mbox{for \ \ all}} \ \ k > N\,.$$
Therefore
$$\Delta_g u (x_k) = c_n \left( u^{{n + 2}\over {n - 2}} (x_k) - u (x_k)
\right)
> \epsilon \ \ \ \ {\mbox{for \ \ all}} \ \ k > N\,,$$
which contradicts (2.2). Therefore $u \le 1$ on $M\,.$ Suppose that there
exists
a point $y_o \in M$ such that $u (y_o) < 1\,.$ Let $v = 1 - u$ and 
$a = \inf_{x \in M} u (x)\,.$ We have $0 < a < 1$. Then $0 \le v \le 1 - a
<
1\,.$ By the Omori-Yau maximum principle, there exists a sequence $\{ y_k
\}_{k
\in {\N}}$ such that
$$\lim_{k
\to \infty} v (y_k) = \sup_{y \in M} v (y)\,, \   \ \ \ \ |\btd v (y_k)|
\le {1\over k}$$ 
and 
$$\Delta_g v (y_k) \le {1\over k}\,. \leqno (2.3)$$ 
As $1 > \sup_{y \in M} v (y) > 0\,,$ there exists a number 
$\delta > 0$ and
$N'
\in {\N}$ such that 
$$\Delta_g v (y_k) =c_n \left( u (y_k) - u^{{n + 2}\over {n - 2}} (y_k)
\right)
> \delta \ \ \ \ {\mbox{for \ \ all}} \ \ k > N'\,,$$ 
which contradicts (2.3). Hence $u \equiv 1$ on $M\,.$ \qedwh
{\bf Theorem 2.4.} \ \ {\it Let $(M, g)$ be a complete non-compact
Riemannian
manifold with scalar curvature equal to $-1\,.$ Assume that the sectional
curvature of $(M, g)$ is bounded from below and the injectivity radius of
$(M,
g)$ is positive. Let $u \in C^\infty (M)$ be a positive solution to
equation
(1.1). If $\inf_{x \in M} u (x) > 0\,,$ then $u \equiv 1$ on
$M$.}\\[0.1in]
{\it Proof.} \ \ By lemma 2.1 we just need to show that $\sup_{x \in M} u
(x) <
\infty\,.$ Let $i_o$ be the injectivity radius of $(M, g)\,.$ By the
assumption
of the theorem, $i_o > 0\,.$ Let $x \in M$ and $\chi \in C^\infty_o (M)$
be a
non-negative function which is zero outside $B_x (i_o)\,.$ As in
\cite{A-M}, we
multiple both sides of (1.1) by $\chi^n u$ and then integrate by parts, we
obtain
\begin{eqnarray*}
(2.5) \ \ \ \ \ \ \ \  & \ & c_n \int_M \chi^n u^{{2n}\over {n - 2}}
\,\,dv_g\\
& \ & 
\le -
\int_M
\chi^n |
\btd u |^2
\,dv_g - n
\int_M u \chi^{n - 1} \left( \btd \chi \cdot \btd u \right) \,dv_g  + c_n
\int_M
\chi^n u^2 \,dv_g\,. \ \ \ \ \ \ \ \ \ \ \ \ 
\end{eqnarray*}
We have 
\begin{eqnarray*} -n\, u \,\chi^{n - 1} \left( \btd \chi \cdot \btd u
\right) & =
& - 2
\left[ {n \over 2} (\chi)^{ {n\over 2} - 1} u \btd \chi \right] \cdot
\left[
(\chi)^{n\over 2} \btd u
\right]\\ & \le & {{n^2}\over 4} u^2 \chi^{n - 2}\, | \btd \chi |^2 \, +
\,\chi^n |
\btd u |^2\,,
\end{eqnarray*} 
where the dot product and the norms are with respect to the Riemannian
metric
$g\,.$ From (2.5) we have
$$c_n \int_M \chi^n u^{{2n}\over {n - 2}} \,\,dv_g \le {{n^2}\over 4}
\int_M u^2
\chi^{n - 2} |
\btd \chi |^2 \,dv_g +  c_n \int_M \chi^n u^2 \,dv_g\,. \leqno (2.6)$$ 
Using Young's inequality we have
$$u^2 \chi^{n - 2} \, | \btd \chi |^2 \le {{(n - 2) \epsilon} \over n}
\,u^{{2n}\over {n - 2}} \,\chi^n + { {2\epsilon^{ - {{n - 2}\over 2}
}}\over n}
\,|
\btd
\chi |^n \leqno (2.7)$$
and
$$\chi^n u^2 \le {{(n - 2) \epsilon} \over n}
\,u^{{2n}\over {n - 2}} \,\chi^n +  {{2\epsilon^{ - {{n - 2}\over 2}
}}\over n}
\, \chi^n\,, \leqno (2.8)$$
where $\epsilon$ is a positive number.  If we choose $\epsilon$ to be
small,
then there exists a positive constant $C
(n)$ depending on $n$ only such that
$$\int_M \chi^n u^{{2n}\over {n - 2}} \,\,dv_g \le  
C (n) \int_M  \left( |\btd \chi |^n + \chi^n \right) dv_g\,. \leqno
(2.9)$$
If we choose $\chi$ to be equal to one on $B_x (i_o/2)$ and zero outside
$B_x
(i_o),$ with $|\btd \chi | \le 10/i_o$, we have
$$\int_{B_x ({{i_o}\over 2}) } u^{{2n}\over {n - 2}} \,\,dv_g \le C (n,
i_o) \
{\mbox{Vol}} \ B_x (i_o)\,.$$
As the sectional curvature is bounded from below and the scalar curvature
is
equal to $-1$, the sectional curvature is bounded from above as well. We
have 
Vol $B_x (i_o)
\le C (i_o),$ where $C (i_o)$ is a positive constant depending on $i_o$
and the
upper bound for sectional curvature only \cite{Gallot-Hulin-Lafontaine}.
Thus
$$\int_{B_x ({{i_o}\over 2}) } u^{{2n}\over {n - 2}} \,\,dv_g \le C$$
and the positive constant $C$ is independent on $x$ and $u\,.$ On $(M, g)$
we
have the Sobolev inequality for compactly supported functions
\cite{Aubin.1}, and
$u$ is a positive smooth function which satisfy 
$$\Delta_g u \ge - c_n u \ \ \ \ {\mbox{on}} \ \ M\,.$$
By using a result of Li-Schoen \cite{Li-Schoen} and the De
Giorgi-Nash-Moser
regularity theory as in \cite{Li},
we conclude that $u$ is bounded from above.{{\hfill
{$\rlap{$\sqcap$}\sqcup$}}}

\vspace{0.5in}

{\bf \Large {\bf 3. \ \  Conformal Transformations}}

\vspace{0.3in}

{\bf Theorem 3.1.} \ \ {\it Let
$(M, g)$ be a complete non-compact Riemannian manifold with constant
negative scalar curvature. Assume that the sectional curvature of
$(M, g)$ is bounded from below and the injectivity radius of $(M, g)$ is
positive. Then any conformal transformation
$F: (M, g) \to (M, g)$ is an isometry.}\\[0.1in] {\it Proof.} \ \ By a
scaling we may assume that the scalar curvature is equal to $-1\,.$ Let
$u = |F'|^{{n - 2}\over 2}\,,$ where $|F'|$ is defined by the equation
$$F^* g \,(y) = |F' (x)|^2 g (x) \ \ \ \ {\mbox{with}} \ \ y = F (x)\,.$$
Then
$u$ satisfies the equation (c.f. \cite{Schoen})
$$\Delta_g u (x) + c_n u (x) = c_n u^{{n + 2}\over {n - 2}} (x) \ \ \ \ 
{\mbox{for \ \ all}} \ \ x \in M\,. \leqno (3.2)$$ 
By the proof of theorem 2.4, $u$ is bounded
from above. On the other hand
$$F^{-1} : (M, g) \to (M, g)$$ 
is also a conformal diffeomorphism and 
$$(F^{-1})^* g (x) = {1\over {|F' (x)|^2 (y)}}\, g (y)  \ \ \ \
{\mbox{with}} \
\  x = F^{-1} (y)\,.$$ Thus $1/u$ also satisfies the equation
$$\Delta_g \left({1\over {u (y)}}\right) + c_n \left({1\over {u
(y)}}\right) =
c_n
\left({1\over {u (y)}}\right)^{{n + 2}\over {n - 2}}
\ \
\
\ {\mbox{for \ \ all}}
\
\ y \in M\,.$$ By the proof of theorem 2.4, $1/u$ is bounded from above.
Therefore 
$u$ is bounded away from zero. By lemma 2.1 we have $u \equiv 1$ on $M$.
Therefore
$F$ is an isometry.\qedwh 
{\bf Corollary 3.3.} \ \ {\it Let $(M, g)$ be a complete
Riemannian  manifold with non-positive sectional curvature and constant
negative
scalar curvature. Then any conformal transformation
$F: (M, g) \to (M, g)$ is an isometry.}\\[0.1in] 
{\it Proof.} \ \ By a scaling we may assume that the scalar curvature is
equal
to $-1\,.$ Let $(\tilde M\,, \tilde g)$ be the universal covering of $M$
with the pull-back metric $\tilde g\,.$ As the scalar curvature of $\tilde
g$
is equal to $-1$ and the sectional curvature is non-positive, therefore
the
sectional curvature is bounded from below. Furthermore, by the
Cartan-Hadamard 
theorem \cite{Ballmann-Gromov-S}, $(\tilde M\,, \tilde g)$ has infinite
injectivity radius. By pulling back the conformal transformation $F$ to a
conformal transformation
$\tilde F :
\tilde M \to \tilde M$ and apply theorem 3.1, $\tilde F$ is an isometry.
Therefore $F$ is also an isometry.\qed
It follows that a conformal transformation of a locally symmetric space of
non-compact type is an isometry. 
\vspace{0.5in}

{\bf \Large {\bf 4. \ \  Complete Solutions}}

\vspace{0.3in}

We remark that without the assumption $\inf_{x \in M} u (x) >
0$, equation (1.1) may have more than one positive solution. 
For example, let $n \ge 3$ be an integer and 
$\delta_{ij}$ be the standard Euclidean metric and 
$$h_{ij} (x) = {{4 n (n - 1) }\over { (1 - |x|^2)^2 }}\, \delta_{ij} \ \ \
\
{\mbox{for}} \ \ x \in B^n \leqno (4.1)$$  be the
(scaled) Poinc\'are metric on the open unit ball ${\bf B}^n$ of
${\R}^n\,,$ 
scaled by the positive number $n (n - 1)$ so that the scalar curvature of
$h$ is $-1\,.$ The (scaled) hyperbolic space
${\Hb}^n = ({\bf B}^n, h)$ is conformal to the standard Euclidean metric
$\delta_{ij}$ on
${\bf B}^n\,.$  
Let $v$ be a positive smooth function and $g_c = v^{4\over {n -
2}} g\,.$ For any smooth function $u$ we have
$$(\Delta_{g} - c_n R_g) (uv) = v^{{n + 2}\over {n - 2}} (\Delta_{g_c} -
c_n
R_{g_c} ) u \leqno (4.2)$$
(see, for example, \cite{Schoen-Yau-Book}). Here $R_g$ and $R_{g_c}$ are
the
scalar curvature of $g$ and $g_c$ respectively. It follows that if $u$ is
a
positive solution to the equation
$$(\Delta_{g_c} - c_n R_{g_c} ) u = c_n u^{{n + 2}\over {n - 2}}\,,$$
then $uv$ satisfies the equation
$$(\Delta_{g} - c_n R_g) (uv) = c_n (uv)^{{n + 2}\over {n - 2}}\,.$$
In particular if $g_{ij} = \delta_{ij}$ and
$$v^{4\over {n - 2}} (x) =  {{4 n (n - 1) }\over { (1 - |x|^2)^2 }}
 \ \ \ \ {\mbox{for}} \ \ x \in B^n\,, \leqno (4.3)$$
and $u$ is a solution to equation (1.1) on the (scaled) hyperbolic space,
then
$$\Delta_o (uv) = c_n (uv)^{{n + 2}\over {n - 2}} \ \ \ \ {\mbox{on}} \ \
B^n\,,$$ where $\Delta_o$ is the Laplacian of the Euclidean space.
Therefore any
positive solution $w$ to the equation 
$$\Delta_o w = c_n w^{{n + 2}\over {n - 2}} \leqno (4.4)$$
also provides a positive solution to equation (1.1) on the (scaled)
hyperbolic
space. Let $(r, \Theta)$ be the polar coordinates on
${\R}^n\,,$ where $r \ge 0$ and $\Theta \in S^{n - 1}\,.$ For a function
$w$
that depends on $r$ only, we have
$$\Delta_o w (r) = w'' (r) + {{n - 1}\over r} w' (r)\ \ \ \ {\mbox{for}} \
\ r >
0\,.$$ 
A calculation shows that equation (4.4) has positive solutions of
the form
$$w_b (x) = {{ \left(4 b^2 n (n - 1) \right)^{{n - 2}\over 4} }\over
{\left(b^2 -
|x|^2\right)^{{n - 2}\over 2} }} \ \ \ \ {\mbox{for}} \ \ \ x \in B^n\,,
\leqno (4.5)$$ 
where $b \ge 1$ is a constant. Thus for any number $b \ge 1\,,$ 
$$u_b (x) = w (x) v^{-1} (x)  = \left( {{b \,(1 - |x|^2)}\over { b^2 -
|x|^2}}
\right)^{{n - 2}\over 2}  \ \ \ \ {\mbox{for}} \ \ \ x \in B^n \leqno
(4.6)$$
is a positive solution to
equation (1.1) on the (scaled) hyperbolic space. We note that if $b >
1\,,$
then $u_b < 1$ and $\lim_{|x|
\to 1} u (x) = 0\,,$ and if
$b = 1\,,$ then $u_b
 \equiv 1$ on $B^n\,.$ By the proof of theorem 2.4, any
positive solution $u$ to equation (1.1) on the hyperbolic space is bounded
from
above. Using the maximum principle as in lemma 2.1 we have $u \le 1$
on the hyperbolic space.\bk
A positive solution to equation (1.1) is called {\it complete} if the
conformal metric $g_c = u^{4\over {n - 2}} g$ is a complete Riemannian
metric
on $M\,.$ In this section we discuss whether complete positive
solutions to the equation (1.1) are unique. We note that the solutions
given in
(4.6) on the (scaled) hyperbolic space are not complete  unless $b = 1$,
that
is, $u \equiv 1\,.$\\[0.2in]      
{\bf Proposition 4.7.} \ \ {\it  Let
$(M, g)$ be a complete non-compact Riemannian manifold with scalar
curvature
equal to $-1\,.$ Assume that the sectional curvature of $(M, g)$ is
bounded from
below and the injectivity radius of $(M, g)$ is positive. Let $u$ be a
positive
smooth solution to the equation (1.1). If the conformal metric
$g_c = u^{4\over {n - 2}} g$ is a complete Riemannian metric with
sectional
curvature bounded from below and with positive injectivity radius, then $u
\equiv 1$ on $M\,.$}\\[0.1in] {\it Proof.} \ \ The Riemannian manifold
$(M,
g_c)$ has constant scalar curvature equal to $-1\,.$ Then $v = 1/u$
satisfies
the equation
$$\Delta_{g_c} v + c_n v = c_n v^{{n + 2}\over {n - 2}} \ \ \ \
{\mbox{on}} \
\ M\,. \leqno (4.8)$$ 
The proof of theorem 2.4 shows that $v$ is bounded from
above, that is, $u$ is bounded away from zero. Lemma 2.1 implies that $u
\equiv 1$ on $M\,.$\qed 
Let us consider positive radial solutions to
equation (1.1). Let
$(r\,,
\Theta)$ be polar coordinates on ${\R}^n\,,$ where $r \in [0, \infty)$ and
$\Theta \in S^{n - 1}\,.$ Let 
$$g = dr^2 + f^2 (r) \,d\Theta^2$$ 
be a Riemannian
metric on
$M\,,$ where
$f \in C^\infty [0, \infty)$ is a positive function with $f (0) = 0$ and
$f' (0) = 1\,.$ We note that $g$ is a complete Riemannian metric on
${\R}^n\,.$ The hyperbolic space is corresponding to
$f (r) =
\sinh r$ for all $r \ge 0\,.$\\[0.2in]
{\bf Theorem 4.9.} \ \ {\it Let $u$ be a positive radial $C^2$-solution to
equation (1.1) on $({\R}^n\,, \ g = dr^2 + f^2 (r) \,d\Theta^2)\,,$ $n \ge
3$.
Assume that  
$\lim_{r \to \infty} f (r) = \infty$ and there
exist positive constants
$R_o$ and $C_o$ such that 
$|f' (r)/ f (r)| \le C_o$ for all
$r > R_o\,.$ If the conformal metric
$g_c = u^{4\over {n - 2}} g$ is complete, then
$u \equiv 1\,.$}\\[0.1in] 
{\it Proof.} \ \ For the Riemannian metric $g  = dr^2 + f^2 (r)
\,d\Theta^2$ we
have
$$\Delta_g = {{\partial^2}\over {\partial r^2}} + {{(n -1) f' (r)}\over {f
(r)}}
{{\partial}\over {\partial r}} + {1\over {f^2 (r)}} \Delta_\Theta\,,
\leqno
(4.10)$$ 
where $\Delta_\Theta$ is the Laplacian for the standard unit sphere in
${\R}^n\,,$ see, for example, \cite{Leung.3}. If $u = u (r)\,,$ then
equation
(1.1) becomes
$$u'' (r) + {{(n -1) f' (r)}\over {f (r)}} u' (r) = c_n (u^{{n + 2}\over
{n - 2}}
(r) - u (r))\,. \leqno (4.11)$$ 
We consider the case that $u (0) = 1$ first. If there exists $r' > 0$ such
that $u (r') \ge  1$ and $u' (r') > 0\,,$ then $u' (r) > 0$ for all $r >
r'\,.$
Otherwise there exists a point $r'' > 0$ such that $u (r'') > 1\,,$ $u'
(r'')
= 0$ and
$u'' (r'') \le 0$, but equation (4.11) shows that this is not
possible. Similarly if there is a point $r_o$ such that $u (r_o) \le 1$
and $u'
(r_o) < 0$, then $u' (r) < 0$ for all $r \ge r_o\,.$ 
Hence if $u \not\equiv 1\,,$ then either
$u (r) > 1$ and $u' (r) > 0$ or $u (r) < 1$ and $u' (r) < 0$ for all $r >
r_o\,.$ Consider the first case. From equation (4.11) we have
$${{(f^{n - 1} u')'}\over {f^{n - 1}}}  = c_n (u^{{n + 2}\over {n - 2}} -
u)\,.
\leqno (4.12)$$ 
Let $\chi \in C^\infty_o ([0, \infty))$ be a cut-off function.
Multipling both sides of equation (4.12) by
$\chi^n u$ and then using integration by parts we obtain
$$c_n \int_0^\infty \chi^n u^2 \,dr -  \int_0^\infty (f^{n - 1} u')
\left({{\chi^n
u}\over {f^{n - 1}}} \right)' \,dr = c_n \int_0^\infty \chi^n u^{{2n}\over
{n -
2}}
\,\,dr\,. \leqno (4.13)$$
We have
$$- (f^n u') \left({{\chi^n u}\over {f^n}} \right)' = -n \chi^{n - 1} u
\chi'
u' - \chi^n |u'|^2 + (n -1) \chi^n u u' {{f'}\over {f}}\,. \leqno (4.14)$$
Applying the Cauchy inequality we have
$$-n \chi^{n - 1} u \chi' u' = - 2 \left( {n\over {\sqrt 2}} \chi^{
{n\over 2}
- 1} u \chi' \right) \left( {1\over {\sqrt 2}} \chi^{n\over 2} u' \right)
\le {{n^2}\over 2} \chi^{n - 2} u^2 |\chi'|^2 + {1\over 2} \chi^n
|u'|^2\,,$$
$$(n -1) \chi^n u u' {{f'}\over {f}} = 2 \left( {{n - 1}\over {\sqrt 2}}
\chi^{n\over 2} u {{f'}\over {f}} \right) \left({ {(n - 1)^2}\over {\sqrt
2}}
\chi^{n\over 2} u'
\right) \le {{(n - 1)^2}\over 2} \chi^n \left({{f'}\over {f}}\right)^2 u^2
+
{1\over 2}
\chi^n |u'|^2\,.$$
Together with (4.13) we obtain 
$$c_n \int_0^\infty \chi^n u^2 \,dr + {{(n - 1)^2}\over 2} \int_0^\infty
\left({{f'}\over {f}}\right)^2 \chi^n u^2 \,dr 
+ {{n^2}\over 2} \int_0^\infty \chi^{n - 2} u^2 |\chi'|^2 \ge 
c_n \int_0^\infty \chi^n u^{{2n}\over {n - 2}} \,dr\,.$$
Apply Young's inequality as in (2.7) and (2.8) and using the bound $|f'/f|
\le
C_o$ to have
$$\int_0^\infty \chi^n u^{{2n}\over {n - 2}} \,dr \le C' \int_0^\infty
(|\chi'|^n
+ \chi^n) \,dr\,, \leqno (4.15)$$
where $C'$ is a positive constant. Let $\chi \equiv 0 $ on $[0, R]
\cup [R+ 3,
\infty)$ with
$R > R_o$ and
$\chi \equiv 1$ on $[R+ 1, R+ 2]\,,$ $\chi \ge 0$ on $[0, \infty)$ and
$|\chi'|
\le {1\over 2}\,.$ From (4.15) we have
$$\int_{R + 1}^{R + 2} u^{{2n}\over {n - 2}} \,dr \le C''$$
for all $R > R_o\,,$ where $C''$ is a constant independent on $R\,.$
Therefore
$u$ is bounded from above. As $u' (r) > 0$ for all $r > 0$ and $u$ is
bounded
from above, we can find a sequence $\{r_k\}$ and a positive constant
$\epsilon > 0$ such that
$u (r_k) > 1 +
\epsilon$ for all $k \in {\N}\,,$ $\lim_{k \to \infty} u' (r_k) = 0$ and
$u'' (r_k) \le 0\,.$ But this contradicts equation (4.11). Thus we must
have $u
(r) < 1$ and $u' (r) < 0$ for all $r > r_o\,.$ Assume that $u (r) \ge c$
for all
$r > \infty\,,$ where $c \in (0, 1)$ is a constant. Then we can find a
sequence
$\{r'_k\}$ and a positive constant
$\delta > 0$ such that
$u (r'_k) > 1 - \delta$ for all $k \in {\N}\,,$ $\lim_{k \to \infty} u'
(r_k)
= 0$ and
$u'' (r_k) \ge 0\,.$ But this contradicts equation (4.11). Therefore
$\lim_{r \to
\infty} u (r) = 0\,.$ There exist positive constants $r' > R_o$ and $C >
0$ such that for
$r \ge r'$ we have
$$(f^{n - 1} u')' (r) = c_n f^{n - 1} (r) (u^{{n + 2}\over {n - 2}} (r) -
u (r))
\le - C f^{n - 1} (r) u (r)\,. \leqno (4.16)$$ Integrating from $r'$ to $r
> R'$
we have
$$f^{n - 1} (r) u' (r) \le f^{n - 1} (r') u' (r') - C \int_{r'}^r f^{n -
1} (s)
u (s) 
\,ds
\le -C u (r)
\int_{r'}^r f^{n - 1} (s)\,ds\,,$$
as $u' \le 0\,.$ Therefore 
$${{u' (r)}\over {u (r)}} \le - C \,{{\int_{r'}^r f^{n - 1} (s) \,ds}\over
{f^{n
- 1} (r)}}\,. \leqno (4.17)$$ 
Using the bounded $f' /f < C_o$ we have $(f^{n - 1})' \le C_o (n - 1) f^{n
-
1}\,.$ An integration gives
$$f^{n - 1} (r) - f^{n - 1} (r') \le C_o (n - 1) \int^r_{r'} f^{n -1}
(s)\,ds\,.
\leqno (4.18)$$ 
As $\lim_{r\to \infty} f (r) = \infty$, if $r$ is large we have
$${1\over 2} f^{n - 1} (r) \le C_o \int^r_{r'} f^{n - 1} (s)\,ds\,,$$
that is,
$${{\int^r_{r'} f^{n - 1} (s)\,ds}\over {f^{n - 1} (r)}} \ge c' \leqno
(4.19)$$
for all $r$ large and for some positive constant $c'\,.$ The inequality
(4.17)
and (4.19) give
$$u (r) \le \tilde Ce^{-cr} \leqno (4.20)$$
for all $r$ large enough, where $\tilde C$ is a positive constant. Thus
the
conformal metric
$g_c = u^{4\over {n - 2}} g$ cannot be complete. Hence $u \equiv 1$ on
$[0,
\infty).$\bk Assume that $u (0) \not= 1\,.$ If $u$ is a $C^2$-function on
${\R}^n$ which depends on $r$ only, then $u' (0) = 0\,.$ If $u (0) > 1\,,$
then
equation (4.11) shows that there exists an $\epsilon > 0$ such that 
$u'' (r) > 0$ for all $r \in (0, \epsilon)\,.$ Hence $u' (r_o) > 0$ for
$r_o$
small enough and hence $u' (r) > 0$ for all $r > r_o\,.$ Similarly if $u
(0) <
1\,,$ then $u' (r) < 0$ for all $r > r_o'\,.$ In either cases we can
obtain
contradiction as above.\qed
We have the following result for non-radial solutions.\\[0.2in]
{\bf Theorem 4.21.} \ \ {\it For an integer $n \ge 4$ let $g = dr^2 + f^2
(r)\,
\,d\Theta^2$ be a Riemannian metric on ${\R}^n\,.$ Assume that $\lim_{r
\to
\infty} f (r) =
\infty$ and there exist positive constants
$R_o$ and $C_o$ such that 
$|f' (r)/ f (r)| \le C_o$ for all
$r > R_o\,.$ Let $u$ be a positive smooth solution to
equation (1.1) on $({\R}^n\,, g)\,.$ If there exist
constants $\delta \in (0, 1)$ and $r_o > 0$ such that $u (r, \Theta) \le 1
-
\delta$ for $r \ge r_o$ and
$\Theta \in S^{n - 1}\,,$ then the conformal metric $g_c = u^{4\over {n -
2}}
g$ is not complete.}\\[0.1in]
{\it Proof.} \ \ In polar coordinates, using (4.10), equation (1.1) is
given by 
$${{\partial^2 u}\over {\partial r^2}} + {{(n -1) f' (r)}\over {f (r)}}
{{\partial u}\over {\partial r}} + {1\over {f^2 (r)}} \Delta_\Theta u
\,=\, c_n
(u^{{n+ 2}\over {n - 2}} - u)\,, \leqno (4.22)$$ 
where $\Delta_\Theta$ is the Laplacian for $S^{n - 1}$ with the standard
metric. For $r \ge r_o$, we have $u (r, \Theta) \le 1 - \delta$ for
some constant $\delta \in (0, 1)$. Therefore 
$$c_n (u^{{n+ 2}\over {n - 2}} - u) \le c_n [ (1 - \delta)^{4\over {n -
2}}  -
1] u = -\epsilon \,u\,, \leqno (4.23)$$
where 
$$\epsilon = c_n [1 - (1 - \delta)^{4\over {n - 2}}] > 0\,.$$
For a fixed number $r \ge r_o\,,$ we integrate equation (4.22) over $S^{n
- 1}$
with respect to the standard measure on $S^{n - 1}$ and use Green's
formula and
(4.23) to obtain
$${{d^2}\over {dr^2}} \left( \int_{S^{n - 1}} u \,d\Theta \right) + {{(n
-1) f'
(r)}\over {f (r)}} {d\over {dr}} \left( \int_{S^{n - 1}} u  \,d\Theta
\right)
\le - \epsilon \int_{S^{n - 1}} u  \,d\Theta\,. \leqno (4.24)$$
For $r \ge r_o\,,$ let 
$$U (r) = \int_{S^{n - 1}} u (r, \Theta) \,d\Theta\,.$$
Then (4.24) can be written as
$$U'' (r) + {{(n -1) f' (r)}\over {f (r)}} U' (r) \le - \epsilon \,U (r) \
\ \ \
{\mbox{for}} \ \ r \ge r_o\,.$$
Therefore
$$\left(f^{n - 1} (r) U' (r)\right)' \le - \epsilon \,f^{n - 1} U (r)  \ \
\ \
{\mbox{for}}
\
\ r
\ge r_o\,.$$ 
The argument in the proof of theorem 4.7 from (4.16) to (4.20) 
shows that there exist positive constants $C$ and $c$ such that 
$$U (r) = \int_{S^{n - 1}} u \,d\Theta \le C e^{- cr}  \ \ \ \
{\mbox{for}}
\
\ r \ge r_o\,. \leqno (4.25)$$ 
Assume that $n \ge 4\,.$ We have 
$$\int_{S^{n - 1}} u^{2\over {n - 2}}  \,d\Theta \le \left( \int_{S^{n -
1}} u 
\,d\Theta \right)^{2\over {n - 2}} \omega_n^{{n - 4}\over {n - 2}}\,,
\leqno
(4.26)$$ 
where $\omega_n$ is the volume of the unit sphere in ${\R}^n\,.$ Using
(4.25)
and (4.26) we have
$$\int_{r_o}^\infty \int_{S^{n - 1}} u^{2\over {n - 2}}  \,d\Theta dr <
\infty\,. \leqno (4.27)$$
We claim that there exists $\Theta_o \in S^{n - 1}$ such that 
$$\int_{r_o}^\infty u^{2\over {n - 2}} (r, \Theta_o) dr < \infty\,.$$
This means that the conformal metric $u^{4\over {n - 2}} (dr^2 + f^2 (r)
\,d\Theta^2)$ is not complete. To prove the claim, by (4.27) and Fubini's
theorem, there exists an positive integer $C'$  such that for any
positive integer
$N$ we have 
$$\int_{r_o}^N \int_{S^{n - 1}} u^{2\over {n - 2}} (r, \Theta) \,d \Theta
dr =
\int_{S^{n - 1}} \left( \int_{r_o}^N u^{2\over {n - 2}} (r, \Theta) \,dr
\right)
d\Theta \le C'\,.$$ 
For each integer $N > r_o$, there exists a point $\Theta_N \in S^{n
- 1}$ such that 
$$\int_{r_o}^N u^{2\over {n - 2}} (r, \Theta_N) dr \le {C'\over
{\omega_n}} +
1\,.
\leqno (4.28)$$ 
A subsequence $\{ \Theta_{N_i} \}_{i \in \N}$ converges to a point
$\Theta_o \in S^{n - 1}\,.$ If 
$$\int_{r_o}^\infty u^{2\over {n - 2}} (r, \Theta_o) dr = \infty\,,$$
then there exists a positive integer $N'$ such that 
$$\int_{r_o}^{N'} u^{2\over {n - 2}} (r, \Theta_o) dr > {C'\over
{\omega_n}}
+ 2\,. \leqno (4.29)$$
As the function 
$$\Theta \to \int_{r_o}^{N'} u^{2\over {n - 2}} (r, \Theta) dr$$
is continuous, therefore in a neighborhood of $\Theta_o$ we have
$$\int_{r_o}^{N'} u^{2\over {n - 2}} (r, \Theta) dr > {C'\over {\omega_n}}
+
{3\over 2}\,. \leqno (4.30)$$
As $\lim_{i \to \infty} \Theta_{N_i} = \Theta_o$ 
and for all $i$ such that $N_i > N'\,,$ $\Theta_{N_i}$ satisfies (4.28)
and
hence
$$\int_{r_o}^{N'} u^{2\over {n - 2}} (r, \Theta_{N_i}) dr \le {C'\over
{\omega_n}} + 1\,.$$
This contradicts with (4.30). The proof of the claim is completed.\qedwh
{\bf Theorem 4.31.} \ \ {\it For an integer $n \ge 3\,,$ let $g = dr^2 +
f^2
(r)\, \,d\Theta^2$ be a Riemannian metric on ${\R}^n\,.$ Assume that
$\lim_{r \to
\infty} f (r) =
\infty$ and there exist positive constants
$R_o$ and $C_o$ such that 
$|f' (r)/ f (r)| \le C_o$ and $|f'' (r) / f (r)| < C_o$ for all
$r > R_o\,.$ Let $u$ be a positive smooth solution to equation (1.1) on
$({\R}^n\,, g)\,.$ If $u \le 1$ on ${\R}^n$ and $u (0) < 1\,,$ then
there exist positive constants
$a$ and $C$ such that $u (r, \Theta) \le 1 - a e^{-C r}$ for all $r \ge
0\,.$}\\[0.1in]
{\it Proof.} \ \ Let $\phi \in C^\infty ({\R}^n)$ be a positive function
such
that $u \le 1 - \phi$ on $B_o (\epsilon)$ and $\phi (r\,, \Theta) = a
e^{-C
r}$ for $r \ge \epsilon$ and $\Theta \in S^{n - 1}\,,$ where $\epsilon$ is
a
small positive constant and
$a$ and $C$ are positive constants to be chosen later. Let $\Phi = u - ( 1
-
\phi)\,.$ If there exists a point $(r_o, \Theta_o) \in {\R}^n$ such that 
$$\Phi (r_o, \Theta_o) = \delta > 0\,,$$
then sup $\Phi \ge \delta\,.$ As the solution $u$ is bounded from above
and the
radial Ricci curvature of $g$ is given by $- (n - 1)f''/ f$, which is
bounded
from below by the assumption on $f$, the maximum principle
(c.f. \cite{Ratto-Rigoli-Setti}) implies that there exists a sequence
$\{(r_k\,,
\Theta_k) \}_{k \in \N}$ such that 
$$\lim_{k
\to \infty} \Phi (r_k\,,
\Theta_k) = \sup \Phi \ge \delta\,, \   \ \ \ \ |\btd \Phi (r_k\,,
\Theta_k)|
\le {1\over k}$$  and 
$$\Delta_g \Phi (r_k\,,
\Theta_k) \le {1\over k}\,.$$
If the sequence is bounded, then we may assume that 
$$\lim_{k \to \infty} (r_k\,, \Theta_k) = (\tilde r\,, \tilde \Theta)\,.$$ 
At
$(\tilde r\,, \tilde \Theta)$ we have
$$\Phi (\tilde r\,, \tilde \Theta) = \sup \Phi \ge \delta \ \ \  \
{\mbox{and}}
\ \ \ \ \Delta_g \Phi (\tilde r\,, \tilde \Theta) \le 0\,.$$
The function $s \mapsto s^{{n + 2}\over {n - 2}} - s$ is increasing on
$[\{(n -
2)/(n + 2)\}^{{n - 2}\over 4}\,, \infty)\,.$ If we choose $a$ to be small,
then
$$u (\tilde r\,, \tilde \Theta) \ge 1 - a e^{-C \tilde r} \ge
\left({{n - 2}\over {n + 2}}\right)^{{n - 2}\over 4}\,. \leqno (4.32)$$
Using (4.10) to compute $\Delta_g \phi$ and (4.32), we have
\begin{eqnarray*}
(4.33) \ \ \ \ \ \ \ \ \ \ \ \ \ \ \ \ \ \ \ \Delta_g \Phi (\tilde r\,,
\tilde
\Theta) & = & c_n (u^{{n + 2}\over {n - 2}} - u) + \Delta_g \phi\\
& \ge & c_n \left[ (1 - a e^{-C \tilde r})^{{n + 2}\over {n - 2}} - (1  -
a
e^{-C \tilde r})\right]\\ & \ & \ \ \ \ \ \  + \,C^2 a  e^{-C \tilde r} -
(n - 1)
{{f' (\tilde r)}\over {f (\tilde r)}} C a e^{-C \tilde r}\,. \ \ \ \ \ \ \
\ \ \
\ \ \ \ 
\end{eqnarray*}
As $f'/f$ is bounded outside $B_o (\epsilon)$, if we choose $C$ to be
large
enough, then there exists a positive constant
$c$ that depends on $C$ and $n$ only, such that 
$$c_n \left[ (1 - a e^{-C \tilde r})^{{n + 2}\over {n - 2}} - (1 - a e^{-C
\tilde r})\right] + C^2 a  e^{-C \tilde r} - (n - 1) {{f' (\tilde r)}\over
{f
(\tilde r)}} C a e^{-C \tilde r} \ge c\,. \leqno (4.34)$$ 
This contradicts that $\Delta_g \Phi
(\tilde r\,, \tilde \Theta) \le 0\,.$ Therefore we may assume that $r_k
\to
\infty$ as $k \to \infty\,.$ Then 
$$\lim_{k
\to \infty} \Phi (r_k\,,
\Theta_k) \ge \delta\,,$$
which implies that $u (r_k\,, \Theta_k) \ge 1 + \delta - a e^{-C r_k} \ge
1 +
\delta/2 > 1$ for $k$ large. This contradicts with the assumption that $u
\le 1$
on
${\R}^n\,.$\qedwh
{\bf Corollary 4.35.} \ \ {\it Let $u$ be a positive solution to equation
(1.1)
on the (scaled) hyperbolic space. If $u \not\equiv 1\,,$ then $u < 1$ on
$B^n\,.$}\\[0.1in]
{\it Proof.} \ \ We know that $u \le 1\,.$ Assume that there is a point
$x_o \in
B^n$ such that
$u (x_o) < 1\,.$ Using an isometry we may assume that $x_o = 0\,.$ The
(scaled) 
hyperbolic metric (4.1) is corresponding to a metric $g = dr^2 + f^2 (r)
\,d\Theta^2$ with 
$$f (r) = \sqrt{n(n - 1)} \sinh {r\over { \sqrt{n(n - 1)} }}\,.$$ 
We may apply theorem 4.31 to complete the proof.\qed

\pagebreak

\end{document}